\documentclass[final,5p,times,twocolumn]{elsarticle}

\usepackage{amssymb}
\usepackage{amsfonts}
\usepackage{latexsym, amssymb,amsmath, amsbsy,amsopn, amstext}
\usepackage{amsmath, amsbsy}
\usepackage{hyperref}
\usepackage{float}
\usepackage{ifpdf}
\usepackage{xcolor}
\usepackage{CJK}
\usepackage{tikz}
\usetikzlibrary{arrows,shapes,snakes,shadows,positioning,automata,patterns}
\usetikzlibrary{trees,decorations.pathmorphing,decorations.markings}

\def\0{{\bf 0}}

\newtheorem{thm}{Theorem}[section]
\newtheorem{lem}[thm]{Lemma}

\newtheorem{rem}[thm]{Remark}

\journal{Systems \& Control letters}

\begin{document}

\begin{frontmatter}



\title{Primal-Dual   Algorithm   for  Distributed Constrained   Optimization \tnoteref{label0}}


 \author{Jinlong ~Lei, Han-Fu Chen, and Hai-Tao Fang  }
 \address{The Key Laboratory  of Systems and Control, Academy of
Mathematics and Systems Science, Chinese Academy of Sciences}
 \ead{leijinlong11@mails.ucas.ac.cn, hfchen@iss.ac.cn, htfang@iss.ac.cn. }
 \begin{abstract}
The paper   studies a distributed constrained optimization  problem, where multiple agents connected in a network
collectively  minimize  the sum of individual objective functions  subject to a global constraint being an intersection of the local constraint sets assigned to  the agents. Based on the augmented Lagrange method, a distributed
 primal-dual algorithm
with  a projection operation included   is proposed  to solve the problem.
It is shown  that with   appropriately  chosen  constant  step size, the local estimates derived at all agents
 asymptotically  reach a consensus at an optimal solution.
In addition,  the value of the cost function at the time-averaged estimate   converges  with rate $O(\frac{1}{k})$
 to the optimal value for the   unconstrained problem. By these properties the proposed
primal-dual algorithm is distinguished from the  existing  algorithms for  distributed constrained optimization.
The theoretical analysis is justified by  numerical simulations.\end{abstract}

\tnotetext[label0]{This
work was supported by the NSFC under Grants 61273193, 61120106011, 61134013, 61174143,
the 973 program of China under grant No.2014CB845301, and
the National Center for Mathematics and Interdisciplinary Sciences, Chinese
Academy of Sciences.}
\end{frontmatter}


\section{Introduction}
 Distributed computation and  estimation   recently   have received much research attention, e.g.,   consensus problems  \cite{murray0, ren},   distributed   estimation  \cite{ Moura1}, sensor localization \cite{localization1},  and  distributed control \cite{control1,control2}.
In particular,  distributed optimization problems have been extensively investigated in  \cite{optimization0}-\cite{Yi}, among which  the distributed subgradient or gradient algorithms  \cite{optimization0}-\cite{Moura} belong to  the  primal domain methods while  \cite{ nedic1}-\cite{Yi} belong to  the primal-dual domain methods.

The paper considers a distributed constrained optimization  problem, where $n$ agents connected in a network collectively minimize the sum of local objective functions $f(x)=\sum_{i=1}^n f_i(x)$ subject to a global constraint   $\Omega_o=\bigcap\limits_{i=1}^n \Omega_i$, where
  $\Omega_i$ is a convex set and $f_i(x)$ is a convex function  in $\Omega_i$. Besides,  $f_i(x)$ and $\Omega_i$  
are  the local data known  to  agent $i$ and  cannot be shared with other agents.
This problem is equivalent to a convex optimization problem with single linear  coupling  constraint and a
convex set constraint. 

 The main contribution of the paper is to   propose a distributed
  primal-dual algorithm with constant step size to solve the constrained optimization problem over the multi-agent network.
  The algorithm is derived  on the basis of   the gradient algorithm
  for finding  saddle points of  an augmented  Lagrange function \cite{Bertsekas}.  
   In an iteration each agent updates its estimate
only using the local data and the information derived from the  neighboring agents.   
 With   appropriately chosen   constant step size, 
the estimates derived at all agents   are shown to reach a consensus at 
 an  optimal solution.   Besides, it is found that  the value of the cost function at the time-averaged estimate converges
 with rate $O(\frac{1}{k})$  to  the optimal  value for   the unconstrained problem.  

 A  general  constrained convex optimization  problem is studied in \cite{nedic1}, where the constraint sets are assumed to be compact. The  problem in the random case is investigated by  \cite{nedic2} for  non-smooth  objective functions,
 meanwhile,   the convex  sets are assumed to be compact and the global constraint set is  required to have a nonempty interior. Here, we study the problem in the deterministic case for  smooth
objective functions, while  imposing  weaker assumptions on the convex sets.

When there are no constraints, the problem of the paper becomes the one discussed in  \cite{optimization0,Moura,Ling,Elia0,Elia}.   The estimates produced by   the  distributed gradient descent (DGD) algorithm  with constant step size   \cite{optimization0} converge to a neighborhood of the optimal  solution.    In contrast to this,  our algorithm gives the accurate estimate. To solve the distributed optimization problems,  
  some  continuous-time distributed  algorithms are proposed  in \cite{Elia0,Elia}, 
 while here   the  discrete-time distributed  algorithm is  investigated.
 The estimates generated  by    the   fast distributed gradient algorithms  \cite{Moura} and by  EXTRA \cite{Ling} converge to an  optimal solution, but  in \cite{Moura} each cost function  is   assumed to be  convex with  gradients
 being bounded and  Lipschitz continuous, while EXTRA \cite{Ling}   only  deals  with   unconstrained problems. Though it is shown by \cite{DSA} that EXTRA   \cite{Ling} is   also a saddle point method,  the augmented   Lagrange function used in \cite{Ling}   is different from ours.  Besides, the convergence rate $O(\frac{1}{k})$ derived here for the unconstrained   case is a new  result.  The primal-dual algorithm proposed in the paper can be seen as an extension of EXTRA  \cite{Ling} to constrained problems.

 The  rest of the paper is organized as follows:   In Section 2,
 some preliminary information about  graph theory and convex analysis is provided and
  the problem is formulated.
 In Section 3, a distributed   primal-dual    algorithm is proposed for solving the problem, while
its convergence is proved in Section 4.
 Two numerical examples are demonstrated  in Section 5, and   some concluding remarks  are given in Section 6.

\section{Preliminaries and Problem Statement} 
We first provide some information about graph theory, convex functions, and convex sets. Then we formulate the distributed   constrained optimization problem to  be investigated.
\subsection{ Graph Theory }
Consider a network of  $n$  agents.  The communication relationship  among  the  $n$
  agents is described by an  undirected graph
   $\mathcal{G } = \{ \mathcal{V }, \mathcal{E }_{\mathcal{G }},\mathcal{A}_{\mathcal{G }}\}$, where  $\mathcal{V }=\{ 1,\cdots,n\}$ is the node set with  node $i$ representing  agent $i$;
 $\mathcal{E }_{\mathcal{G }}  \subset  \mathcal{V } \times \mathcal{V } $ is the undirected  edge set,
  and  the unordered pair of nodes   $(i,j)\in\mathcal{E }_{\mathcal{G }}$ if and only if  agent  $i$ and  
  agent $j$ can exchange information with each  other;
 $\mathcal{A}_{\mathcal{G }} =[a_{ij} ]\in  \mathbb{R}^{n\times n} $ is the  adjacency matrix of $\mathcal{G } $,
where  $a_{ij} =a_{ji}>0$  if  $(i,j)\in \mathcal{E }_{\mathcal{G }} $, and  $a_{ij}=0$, otherwise.
Denote by $\mathcal{N}_i= \{ j \in \mathcal{V}: (i,j) \in \mathcal{E}_{\mathcal{G}}\}$ the neighboring agents of   agent $i$.  The Laplacian matrix of   graph   $\mathcal{G}$ is defined as
 $\mathcal{L}_{\mathcal{G}}= \mathcal{D}_{\mathcal{G}}-\mathcal{A}_{\mathcal{G}}$,
where $\mathcal{D}_{\mathcal{G}}=diag\{ \sum_{j=1}^n a_{1j}, \cdots, \sum_{j=1}^n a_{nj} \}$.
For a given pair $i,j \in \mathcal{V}$,  if there exists a sequence of distinct nodes
 $i_1, \cdots, i_{p }$ such  that $ (i,i_1)\in \mathcal{E }_{\mathcal{G }} , $ $ (i_1,i_2)\in \mathcal{E }_{\mathcal{G }} $,
  $\cdots$, $(i_{p},j) \in \mathcal{E }_{\mathcal{G }} $,
then  $(i,i_1, \cdots, i_p,j)$ is called the undirected  path between  $i $ and  $j$.  We say that $\mathcal{G } $ is connected
  if there exists  an undirected path  between  any    $i,j \in \mathcal{V}.$

 The following lemma presents  some properties of  the Laplacian matrix $\mathcal{L}  $  corresponding to an undirected   graph $\mathcal{G}  $.

\begin{lem}  \label{lem1}  \cite{graph}
The Laplacian  matrix   $\mathcal{ L}  $   of  an  undirected graph $\mathcal{G} $ has the following properties:

  i)   $\mathcal{L}$   is symmetric and positive semi-definite;

 ii)   $\mathcal{ L}$  has a simple zero eigenvalue with  corresponding eigenvector  equal  to  $\mathbf{1}$,
 and all other eigenvalues are positive if and only if the graph $\mathcal{G} $ is connected,
 where  $\mathbf{1}$  denotes  the vector of compatible dimension with all entries equal to 1.

\end{lem}

\subsection{Gradient, Projection  Operator  and  Normal Cone}
For  a given function $f: \mathbb{R}^m \rightarrow [-\infty, \infty],$
  denote its  domain  as   $ \textrm{dom} (f) \triangleq \{ x \in \mathbb{R}^m:  f(x) < \infty\}.$
Let $f(\cdot)$ be a convex function, and let $x \in \textrm{dom} (f)$.
For  a smooth (differentiable) function $f(\cdot)$,  denote  by $\nabla f(x)$ the gradient of  the function $f(\cdot)$ at point $x$.
Then  \begin{equation} \label{e1}
  f(y) \geq f(x) + \langle  \nabla f(x), y-x  \rangle~~  \forall y \in \textrm{dom} (f),
\end{equation}
  where $ \langle x,y  \rangle $  denotes the inner product of  vectors $x$ and $y.$

 For a  nonempty convex   set $\Omega \subset \mathbb{R}^m$ and a point $x \in \mathbb{R}^m$, we call the point in
$\Omega$ that is closest  to $x$ the projection of $x$ on $\Omega$ and   denote it by $P_{\Omega} \{ x\}$.
If  $\Omega \subset \mathbb{R}^m$ is  closed,
then  $P_{\Omega} \{ x\}$ contains only one element for any $x \in  \mathbb{R}^m.$

Consider  a convex closed set $\Omega \subset \mathbb{R}^m$ and a point $x \in  \Omega$. Define
   the normal cone  to $\Omega$ at $x$ as   $N_{\Omega} \{ x\} \triangleq \{ v \in \mathbb{R}^m: \langle v, y-x \rangle \leq 0  ~~\forall y \in \Omega\}$.  It is shown in  \cite[Lemma 2.38]{opt}  that    the following equation holds for any   $x \in  \Omega$:
 \begin{equation}\label{normalcone}
N_{\Omega} \{ x\} =\{ v \in \mathbb{R}^m: P_{\Omega} \{x+v\}=x \}.
\end{equation}

   A set  $C$ is affine if it contains the lines that pass through  any
 pairs of points $x,y \in C$ with $x \ne y$. Let   $\Omega \subset \mathbb{R}^m$ be a nonempty convex set.
  We say that $x \in \mathbb{R}^m$ is   a   relative interior point of   $\Omega $ if 
 $x \in \Omega$ and  there exists an open sphere $S$ centered at $x$  such  that
$S \cap \textrm{aff}(\Omega) \subset \Omega,$
where   $\textrm{aff}(\Omega)$  is the intersection of all affine sets containing $\Omega$.
  A point   $x \in \mathbb{R}^m$ is called the    interior  point  of  $\Omega  $ if  $x \in \Omega$ and  
  there exists an open sphere $S$ centered at $x$  such  that   $    S \subset \Omega  .$ 
A  pair of vectors  $ x^{*} \in  \Omega  $ and $  z^{*} \in  \Psi$      is  called a saddle point
 of the function $ \phi(x,  z)$ in $ \Omega \times \Psi$ if  \begin{equation}
\phi(x^{*}, z) \leq \phi(x^{*}, z^{*})  \leq \phi(x , z^{*}) ~~ \forall  x \in \Omega,~ ~ \forall  z \in \Psi. \nonumber
\end{equation}
    These  definitions can be found in \cite{Bertsekas}.

\subsection{ Problem Formulation} 
 Consider a network of $n$ agents that  collectively solve the  following constrained optimization   problem
\begin{equation}\label{problem1}
\begin{split}
&  \textrm{minimize} ~~ f(x)=\sum_{i=1}^n  f_i(x) \\
&  \textrm{subject to} ~~ x \in  \Omega_o= \bigcap_{i=1}^n \Omega_i,
\end{split}
\end{equation}
where  $\Omega_i \subset \mathbb{R}^m$
 is a closed convex set, representing the local constraint set of   agent $i$, 
 and $f_i(x): \mathbb{R}^m\rightarrow \mathbb{R}$ is a smooth convex function in $ \Omega_i$,  representing the local
objective function of   agent $i$.
Assume that $f_i$ and $\Omega_i$ are  privately    known to    agent  $i$.
  We assume that there  exists at least one finite solution $x^*$ to the problem \eqref{problem1}.
For   the problem \eqref{problem1},  denote by   $f^{*}= \min_{ x \in     \Omega_o} f(x)$ the optimal value,
and  by $  \Omega_{o}^{*}=\{  x \in   \Omega_o: f(x)=f^{*} \} $  the optimal solution set.

  We  use an undirected graph  $\mathcal{G } = \{ \mathcal{V }, \mathcal{E }_{\mathcal{G }},\mathcal{A}_{\mathcal{G }}\}$ to describe the communication    among   agents.  
Let  $\mathcal{L}$ denote  the   Laplacian  matrix of  the undirected graph $\mathcal{G }$.

    Let us introduce  the following    conditions for the problem.
\begin{description}
  \item[A1]   $ \Omega_o$     has at least  one relative   interior point.
  \item[A2]  The undirected graph $\mathcal{G }$  is connected.
  \item[A3] For any $i \in \mathcal{V}$,   $\nabla  f_i(x)$ is locally Lipschitz continuous on  $\Omega_i$.
\end{description}

  \section{Algorithm Design}

We first  give an equivalent  form    of the problem \eqref{problem1}. Then
define  a  distributed primal-dual    algorithm with constant step size to solve the formulated  problem.

\subsection{An Equivalent  Problem}
\begin{lem} \label{lem2} \cite[Lemma 3]{Liu}
If A2 holds, then  the problem \eqref{problem1} is equivalent to the following  optimization   problem
\begin{equation}\label{problem2}
\begin{split}
& \textrm{minimize}  ~~\widetilde{ f}(X)=\sum_{i=1}^n  f_i(x_i) \\
& \textrm{subject to} ~~ ( \mathcal{L} \otimes \mathbf{I}_m ) X=\mathbf{0}, ~~  X \in \Omega ,
\end{split}
\end{equation}
 where $X= col\{x_1,\cdots, x_n\}\triangleq  (x_1^T, \cdots, x_n^T)^T$,
$\Omega=\prod_{i=1}^n \Omega_i$ denotes the Cartesian product,
$\otimes $ denotes  the Kronecker product,   $\mathbf{I}_m$  denotes the identity matrix of size $m$,
and $\mathbf{0}$ denotes the vector of compatible dimension with all entries equal to 0.
\end{lem}

\begin{rem}\label{r1}
 Lemma \ref{lem2}    implies    that  solving  the      problem \eqref{problem1}
 is equivalent  to solving  the problem  \eqref{problem2} when the underlying graph is  undirected and  connected.
 If  $X^{*}=col\{ x_1^{*}, \cdots, x_n^{*})$ is a  solution to
  the problem \eqref{problem2},  i.e.,  $ f^{*}=\widetilde{ f}(X^{*})$,
   then $ x_i^{*}= x_j^{*}=x^{*}~\forall i, j\in \mathcal{V}$   for some $x^{*} \in \Omega_o$  by $( \mathcal{L} \otimes \mathbf{I}_m ) X^{*}=\mathbf{0}$ and $X^* \in \Omega$.  Thus,    $\widetilde{ f}(X^{*})=\sum_{i=1}^n f_i(x^{*})=f^{*}$, and hence $x^{*}$
      is   an optimal solution to  the  problem   \eqref{problem1}.
 \end{rem}

Define   the Lagrange function  $ \phi(X,\Lambda)= \widetilde{f}(X)+\langle \Lambda, ( \mathcal{L} \otimes \mathbf{I}_m ) X \rangle$, where  $\Lambda \in  \mathbb{R}^{mn}$ is the Lagrange  multiplier vector.   Then the  original problem  \eqref{problem2} can   be rewritten as $ \inf\limits_{X \in \Omega} \sup\limits_{\Lambda \in \mathbb{R}^{ mn} } \phi(X, \Lambda)$, while  the  dual problem is defined as follows 
\begin{equation}\label{pdual}
\sup\limits_{\Lambda \in \mathbb{R}^{ mn} } \inf\limits_{X \in \Omega}  \phi(X, \Lambda).
\end{equation}

\begin{lem}\label{lem3} Assume A1 and A2  hold. 
   Then $ \phi(X, \Lambda)$ has at least one saddle point in $ \Omega \times \mathbb{R}^{mn}$.
 A  pair  $(X^{*}, \Lambda^{*}) \in  \Omega \times \mathbb{R}^{mn}$  is
the primal-dual  solution   to the problems  \eqref{problem2} and \eqref{pdual}
  if and only if    $(X^{*}, \Lambda^{*}) $ is  a saddle point
 of $ \phi(X, \Lambda)$ in $ \Omega \times \mathbb{R}^{mn}$.
\end{lem}
{\bf Proof}:
Since $f_i(\cdot)~ \forall i\in \mathcal{V}$ are continuous and the problem \eqref{problem1} has at least one finite 
solution,   $f^{*}$ is finite.  Moreover,  A1  implies that there exists  a  relative interior $\bar{X} $ of  set
$\Omega$ such that $(\mathcal{L}  \otimes \mathbf{I}_m) \bar{X}=0$. 
Then  by    \cite[Proposition   5.3.3]{Bertsekas} we know that  the primal and dual optimal  values are equal, i.e.,
 \begin{equation}\label{minmax}
\inf_{X \in \Omega} \sup_{\Lambda \in \mathbb{R}^{ mn} } \phi(X, \Lambda)=\sup_{\Lambda \in \mathbb{R}^{ mn} } \inf_{X \in \Omega}  \phi(X, \Lambda),
\end{equation}and there  exists  at least one dual optimal solution.  So,  by \eqref{minmax} we conclude that $ \phi(X, \Lambda)$ has at least one saddle point in $ \Omega \times \mathbb{R}^{mn}$.

Since the  minimax equality \eqref{minmax} holds,   by   \cite[Proposition 3.4.1]{Bertsekas}  we know that
 $X^{*}$ is  the  primal optimal   solution and $\Lambda^{*}$  is the dual optimal   solution if and only if
 $(X^{*}, \Lambda^{*})$ is  a saddle point  of $ \phi(X,\Lambda)$ on  $\Omega \times \mathbb{R}^{mn}$.
  This completes the proof. \hfill $\blacksquare$

\subsection{ Distributed Primal-Dual Algorithm}

Denote  by $x_{i,k} \in \mathbb{R}^m$ the  estimate for the optimal  solution to  the  problem \eqref{problem1} given by  agent $i$ at time $k$,   and by $\lambda_{i,k} \in \mathbb{R}^m$ the corresponding Lagrange  multiplier. They are updated  as follows:
 \begin{equation} \label{fixed}
 \begin{array}{ll}
   x_{i ,k+1} & = P_{\Omega_i} \{  x_{i ,k} - \alpha   \nabla f_i(x_{i,k})   -
 \alpha   \sum_{j=1}^n   a_{ij}  ( \lambda_{i ,k}-\lambda_{ j ,k}) \\ & \quad\qquad\qquad \qquad\qquad  -\alpha \sum_{j=1}^n   a_{ij} (x_{i,k}-x_{j,k})  \}, \\
 \lambda_{i ,k+1}&=\lambda_{i,k}+\alpha \sum_{j=1}^n   a_{ij} (x_{i,k}  -x_{j ,k}) .
\end{array}
\end{equation}

Set  $ X_k = col\{  x_{1,k} , \cdots, x_{n,k}\}, $ $ \Lambda_k = col\{ \lambda_{1 ,k} \cdots, \lambda_{n ,k}\},$ and
  $   \nabla \widetilde{f} (X_k)= col\{   \nabla f_1(x_{1,k})   ,\cdots,   \nabla f_n(x_{n, k})  \}.$ 
  Then  \eqref{fixed} can be written in the  compact form as follows:
   \begin{align}
& X_{k+1}= P_{\Omega} \{  X_k -  \alpha     \nabla \widetilde{f} (X_k)   -\alpha ( \mathcal{L}   \otimes \mathbf{I}_m)  \big( \Lambda_k +  X_k \big)  \}, \label{primal} \\
&\Lambda_{k+1}=\Lambda_k+\alpha   ( \mathcal{L}  \otimes \mathbf{I}_m)  X_k . \label{dual}
\end{align}

 Note that  the algorithm  \eqref{primal} \eqref{dual} actually  is  the  gradient algorithm for finding saddle points of the augmented Lagrange function $  \widetilde{\phi}(X,\Lambda)=   \phi(X,\Lambda) 
+ \frac{1}{2} \langle X, ( \mathcal{L} \otimes \mathbf{I}_m ) X \rangle$ in
$ \Omega \times \mathbb{R}^{mn}$.
By Lemma \ref{lem3} we see that if  the algorithm  \eqref{primal} \eqref{dual}  converges to a saddle point of the augmented
Lagrange function, then it solves the original problem \eqref{problem2}. Convergence properties of the primal-dual
 method have been studied extensively, see, for example,  \cite{saddle1,saddle2}. In general, only 
  a subsequence of the sequence  $(X_k, \Lambda_k)$ converges to a saddle point of the augmented
Lagrange function. To obtain  the convergence of the whole sequence $(X_k, \Lambda_k)$, it is often to
 assume that the  augmented Lagrangian function is strictly convex-concave.  
  However,   for the problem  studied in the paper, the augmented Lagrange function is neither strongly convex in $X \in  \Omega $ nor strongly concave in $\Lambda \in\mathbb{R}^{mn}$. Thus, the  standard analysis  
  of     gradient methods for finding saddle points  is  not  applicable here. Instead, we apply  the Lyapunov  function method      to  analyze  convergence.

\section{Convergence Analysis}
Convergence  results for the proposed  primal-dual  algorithm  are presented in Section 4.1 with  the   proof
given in Sections 4.2 and  4.3.

   \subsection{ Main  Results}

By  A2  from  Lemma \ref{lem1} we know that all eigenvalues of $\mathcal{L}$
are    nonnegative real numbers, and   zero is  a simple eigenvalue. Let  us write
 the eigenvalues of $\mathcal{L}$     in the  non-decreasing
 order as $0=\kappa_1 <  \kappa_2 \leq \cdots \leq \kappa_n$.
 
Set    \begin{equation}\label{wmatrix}
 \mathcal{W}= (\mathbf{I}_{n } - \alpha\mathcal{L} +\alpha^2 \mathcal{L}^2  ) \otimes \mathbf{I}_m.
 \end{equation}
 Let $(X^*,\Lambda^*)$ be a saddle point of $\phi(X,\Lambda)$.
 Define $$V_1(X)= \langle X-X^*, \mathcal{W} (X-X^*)  \rangle,~~
 V_2(\Lambda)= \| \Lambda-\Lambda^*\|^2.$$
 Construct a candidate Lyapunov function  as follows
 \begin{equation}
\label{ }
V(X, \Lambda )=V_1(X)+V_2(\Lambda). \nonumber
\end{equation}

The following theorem shows   that the local estimates derived at all agents
asymptotically reach a consensus at  an  optimal solution  to the problem \eqref{problem1}.

  \begin{thm}\label{thm1}   Assume A1-A3 hold. Let $\{ x_{i,k}\}$ and  $\{ \lambda_{i,k}\}$   be produced by \eqref{fixed}
  with initial values $x_{i,0},~ \lambda_{i,0}$.      Let   $(X^{*}, \Lambda^{*}) $ be  a saddle point
 of $ \phi(X, \Lambda)$ in $ \Omega \times \mathbb{R}^{mn}$.    Assume,  in addition,  that   the constant step  size $\alpha$ satisfies
   $ 0< \alpha \leq \frac{1}{2\kappa_n}  \textrm{~and ~}\alpha < \frac{3 }{2l_r}$,  where
  $ l_r$ is    the local Lipschitz  constant  of  $\nabla \widetilde{f} (X) $
 in  the compact set $\{X \in \Omega:\| X-X^{*}\| \leq r \} $   with $r$ defined by
\begin{equation} \label{radius}
 r =\sqrt{V(X_0, \Lambda_1) /\lambda_{min}( \mathcal{M})} ,
\end{equation}
where  $\lambda_{min}(  \cdot)$ denotes the smallest eigenvalue of a symmetric matrix, 
and $\mathcal{M}=diag\{I_{mn}, \mathcal{W}\}.$   Then 

(i) $   V(X_k, \Lambda_{k+1})    $  monotonously   decreases and converges,

(ii) $ d_k=\sqrt{\| X_k-X^*\|^2 + \| \Lambda_{k+1}-\Lambda^{*}\|^2} \leq r ~~\forall k \geq 0$,  

(iii) $ \lim\limits_{k \rightarrow \infty} x_{i,k}=\lim\limits_{k \rightarrow \infty} x_{j,k}=x^* ~~\forall i,j \in \mathcal{V}$ for some
$x^* \in \Omega_o^*.$  
 
  \end{thm}

\begin{rem}
  The   problem considered    in  \cite{Liu} is in the same form as the problem \eqref{problem1},
 but the local constraint is  a hyper-box or hyper-sphere,  which  is  a special case of   A1.  Unlike the discrete-time  algorithm  \eqref{fixed},  the  continuous-time  distributed algorithm is proposed   in   \cite{Liu}.
     Though the estimates   given by all agents  converge  to the same optimal solution,     some  intermediate  sequence  might be unbounded,
 which  makes the algorithm difficult  to be implemented.
   \end{rem}

 Denote by $\bar{X}_k=\frac{1}{k+1}\sum_{p=0}^k X_p$  the time-averaged  estimate. 
   In what follows,  the   convergence  rate   of the  algorithm \eqref{fixed} 
  for the case where    $\Omega_i =\mathbb{R}^m~~\forall i\in \mathcal{V}$ is shown.

  \begin{thm} \label{thm2} Assume  $\Omega_i =\mathbb{R}^m~\forall i\in \mathcal{V}$,  A2, and   A3 hold.
  Let $\{ x_{i,k}\}$ and  $\{ \lambda_{i,k}\}$   be produced by  the algorithm \eqref{fixed}
 with    initial values $x_{i,0},~ \lambda_{i,0}$.
  Let  $(X^{*},\Lambda^{*}) $ be   a saddle point of $\phi(X,\Lambda) $   in $ \Omega \times \mathbb{R}^{mn}$.
  If $ 0< \alpha \leq \frac{1}{2\kappa_n}  \textrm{~and ~}\alpha < \frac{3 }{2l_r}$, where
  $ l_r$ is  the local Lipschitz  constant  of  $\nabla \widetilde{f} (X) $
 in  the compact set $\{X \in \Omega:\| X-X^{*}\| \leq r \} $
 with $r $ defined by \eqref{radius},   then
 \begin{flalign}
 & \textrm{(i) } ( \mathcal{L}  \otimes \mathbf{I}_m)  \bar{X}_k= \frac{\Lambda_{k+1}-\Lambda_0}{(k+1)  \alpha   }, \label{con}\\
&  \textrm{(ii)}  ~ \widetilde{f}(\bar{X}_k) \leq f^{* } + \frac{1}{2\alpha (k+1)}  \Big(  \|  X_0 -X^{*}\|^2 +\| \Lambda_0\|^2- \|  X_{k+1} - X^{*}\|^2 \nonumber \\&   -\| \Lambda_{k+1}\|^2  \Big)  +  \frac{c_r }{2\alpha(k+1)} \times   \Big(   V(X_0, \Lambda_{1})- V(X_{k+1}, \Lambda_{k+2})      \Big) ,  \label{sum0}  \\
&\textrm{(iii)}   ~\widetilde{f}(\bar{X}_k)  \geq   f^{*}-\frac{ \langle \Lambda_{k+1}-\Lambda_0 , \Lambda^{*} \rangle }{(k+1)  \alpha   } , \label{sum00} \end{flalign}
where   $c_r=   1/\lambda_{min} \big( \mathcal{W} -\frac{\alpha  l_r}{2 } I_{mn} \big) +1   $.
\end{thm}

      \begin{rem}
   Since $ \| \Lambda_{k}-\Lambda^{*}\|  \leq r~\forall k \geq 0$ by Theorem \ref{thm1}(ii),  $\Lambda_{k}$ is uniformly 
bounded in $k$ by a constant.  Then by \eqref{con}  we  see  that  $ ( \mathcal{L}  \otimes \mathbf{I}_m)  \bar{X}_k$   converges to 
$\mathbf{0}$ with   rate $O(\frac{1}{k}). $   Since  Theorem \ref{thm1}(i)  implies that   $  V(X_0, \Lambda_{1}) - V(X_{k+1}, \Lambda_{k+2})  \leq 0~~\forall k \geq 0, $  by \eqref{sum0}\eqref{sum00}   the value of  the cost function $\widetilde{f}(\cdot)$
   at   $\bar{X}_k$ converges  to the  optimal value  with rate $O(\frac{1}{k}). $ \end{rem}

\begin{rem} Note that $ \Lambda_{1}=\Lambda_0+\alpha   ( \mathcal{L}  \otimes \mathbf{I}_m)  X_0.$
Then from Theorems \ref{thm1} and \ref{thm2}  we see that for   small enough $\alpha>0$ depending   on the distance between the initial value and the optimal solution,   and on the structure of the cost functions
in the neighborhood of the optimal solution,  the  estimates  given by all
agents  finally  reach a  consensus at an  optimal solution.
If   $\nabla  f_i(x)$ is globally  Lipschitz continuous in set $\Omega_i$ with constant  $l_c$ for  any $i \in \mathcal{V}$,
then    the results given in Theorems \ref{thm1} and \ref{thm2}    hold as well  for any  $\alpha$ satisfying  $ 0< \alpha \leq \frac{1}{2\kappa_n}  \textrm{~and ~}\alpha < \frac{3 }{2l_c}$ but  independent  of  the initial values.
\end{rem}
 \subsection{Proof of Theorem \ref{thm1}}
 Prior to proving  Theorem  \ref{thm1},  we give a lemma that will be used in the proof.

 \begin{lem}  \cite[Theorem 2.1.5]{nesterov}  \label{lemma5}
 If $f: \mathbb{R}^m  \rightarrow    \mathbb{R}$ is a convex function whose gradient is   globally Lipschitz
 continuous   with constant $l_c$,  then
 $$ \langle  x-y , \nabla f(x) - \nabla f(y)\rangle  \geq  \frac{1}{l_c } \| \nabla f(x)- \nabla f(y) \|^2 ~~\forall x,y \in \mathbb{R}^m.$$
 \end{lem}

\textbf{ Proof  of Theorem \ref{thm1}:}  
   Note that
 \begin{equation}  \label{duct0}
\begin{array}{lll}
&    V(X_{k+1}, \Lambda_{k+2})-V(X_k, \Lambda_{k+1})\\
&=V_1(X_{k+1})-V_1(X_k)+
V_2(\Lambda_{k+2})- V_2(\Lambda_{k+1}) \\& =  \langle  \mathcal{W}  ( X_{k+1}- X_k), X_{k+1}+X_k -2X^{*}  \rangle \\& \quad+ \langle   \Lambda_{k+2}-\Lambda_{k+1} ,  \Lambda_{k+2}+\Lambda_{k+1} -2\Lambda^{*}  \rangle \\&  = -\| \Lambda_{k+2}-\Lambda_{k+1}\|^2  -  \langle X_{k+1}- X_k, \mathcal{W}  ( X_{k+1}- X_k)  \rangle \\&
\quad + 2 \langle  \Lambda_{k+2}-\Lambda_{k+1},  \Lambda_{k+2}-\Lambda^{*} \rangle \\& \quad+ 2  \langle  X_{k+1}-X^*, \mathcal{W} (X_{k+1} -X_k) \rangle    . \end{array}
\end{equation}
We now   estimate  the last two terms on the right hand side of  \eqref{duct0}.

Since $(X^{*},\Lambda^{*}) $ is a saddle point of $\phi(X,\Lambda)$,  by Lemma \ref{lem3} we see  that $X^{*}$ is an optimal solution to the problem \eqref{problem2}, and hence
 \begin{equation}\label{consensus}
( \mathcal{L} \otimes \mathbf{I}_m)  X^{*}=\mathbf{0}.
\end{equation}
Since $ \mathcal{L}$ is symmetric, by  \eqref{dual}  \eqref{consensus}
 we derive
\begin{equation}
\begin{array}{lll}
&\langle  \Lambda_{k+2}-\Lambda_{k+1},  \Lambda_{k+2}-\Lambda^{*}  \rangle
\\& = \langle   \alpha ( \mathcal{L} \otimes \mathbf{I}_m)(X_{k+1} -X^{*}),  \Lambda_{k+2}-\Lambda^{*} \rangle
  \\&= \langle   \alpha ( \mathcal{L} \otimes \mathbf{I}_m)( \Lambda_{k+2}-\Lambda^{*} ), X_{k+1} -X^{*}\rangle   . \nonumber
\end{array}
\end{equation} 
   Thus, 
\begin{equation}\label{inequality01}
\begin{array}{lll}
&  \langle  \Lambda_{k+2}-\Lambda_{k+1},  \Lambda_{k+2}-\Lambda^{*}  \rangle
 +   \langle  X_{k+1}-X^*, \mathcal{W} (X_{k+1} -X_k)  \rangle
 \\& =  \langle    X_{k+1} -X^{*} , \alpha ( \mathcal{L} \otimes \mathbf{I}_m) (  \Lambda_{k+2}-\Lambda^{*}) +  \mathcal{W} ( X_{k+1}-X_k) \rangle  .
\end{array}
\end{equation}

From \eqref{dual} we derive  $\Lambda_{k+2}=\Lambda_k+\alpha   ( \mathcal{L}  \otimes \mathbf{I}_m)  X_k+ \alpha   ( \mathcal{L}  \otimes \mathbf{I}_m)  X_{k+1} $,  and hence by  \eqref{consensus}
 \begin{equation}\label{ex}
\Lambda_{k}-\Lambda^{*}=\Lambda_{k+2}-\Lambda^{*}- \alpha   ( \mathcal{L}  \otimes \mathbf{I}_m)  ( X_{k} -X^{*} )- \alpha   ( \mathcal{L}  \otimes \mathbf{I}_m) ( X_{k+1} -X^{*} ).
\end{equation}
Then   by multiplying  both sides of \eqref{ex} with $( \mathcal{L}  \otimes \mathbf{I}_m) $,  from  the rule of  the Kronecker product $(A \otimes B)(C \otimes D)= (AC ) \otimes ( BD)$   we  obtain
\begin{equation}\label{lamada}
\begin{array}{lll}
&( \mathcal{L}  \otimes \mathbf{I}_m)   ( \Lambda_k  - \Lambda^{*}) =( \mathcal{L}  \otimes \mathbf{I}_m)  ( \Lambda_{k+2}  -\Lambda^{*} ) \\
&- \alpha   ( \mathcal{L}^2  \otimes \mathbf{I}_m)  ( X_k -X^{*} )- \alpha   ( \mathcal{L}^2  \otimes \mathbf{I}_m) ( X_{k+1} -X^{*} ).
\end{array}
\end{equation} 
 
 Set  \begin{equation}\label{z}
Z_{k+1}= X_k-  \alpha      \nabla \widetilde{f} (X_k)    -\alpha ( \mathcal{L}   \otimes \mathbf{I}_m)  \big( \Lambda_k +  X_k \big) -X_{k+1}.
 \end{equation}
Then by \eqref{consensus}  \eqref{lamada}  \eqref{z}   we derive
\begin{equation}\label{recursive}
\begin{array}{lll}
 & X_{k+1}- X^{*} =  X_k - X^{*}- \alpha  \nabla \widetilde{f} (X_k)
  -\alpha ( \mathcal{L}   \otimes \mathbf{I}_m) \big(  X_k-X^{*} \big)  \\&-\alpha ( \mathcal{L}   \otimes \mathbf{I}_m)  \big( \Lambda_k -\Lambda^{*} \big)
  -\alpha ( \mathcal{L}  \otimes \mathbf{I}_m)  \Lambda^{*}   - Z_{k+1} \\
  &=X_k -X^{*}  -  \alpha  \big(\nabla  \widetilde{f} (X_k) - \nabla \widetilde{f} (X^{*})  \big)
   -\alpha ( \mathcal{L}   \otimes \mathbf{I}_m)  \big(    X_k-X^{*} \big)\\&
 - \alpha(\nabla \widetilde{f} (X^{*})+( \mathcal{L}  \otimes \mathbf{I}_m) \Lambda^{*} )   -Z_{k+1} \\
 &-  \alpha \Big ( ( \mathcal{L}  \otimes \mathbf{I}_m)  \big(  \Lambda_{k+2} - \Lambda^{*} \big)
- \alpha   ( \mathcal{L}^2  \otimes \mathbf{I}_m)  ( X_k -X^{*} )
\\& - \alpha ( \mathcal{L}^2  \otimes \mathbf{I}_m)  ( X_{k+1} -X^{*})  \Big)
\\ & = \big( (\mathbf{I}_{n } - \alpha\mathcal{L} +\alpha^2 \mathcal{L}^2  ) \otimes \mathbf{I}_m \big)  (X_k-X^{*})-Z_{k+1}
\\&-\alpha ( \mathcal{L}  \otimes \mathbf{I}_m) \big(  \Lambda_{k+2} - \Lambda^{*} \big)-\alpha \big(\nabla \widetilde{f} (X_k) - \nabla  \widetilde{f} (X^{*})  \big)\\&
 - \alpha(\nabla \widetilde{f} (X^{*})+(\mathcal{L} \otimes \mathbf{I}_m) \Lambda^{*} )   +  \alpha^2
 ( \mathcal{L}^2 \otimes \mathbf{I}_m) ( X_{k+1} -X^{*} ) .
\end{array}
\end{equation}
Moving the last term at  the
right-hand side of \eqref{recursive} to the left and subtracting  $\mathcal{W}(  X_{k+1}- X^{*})$  from  both  sides of \eqref{recursive} we derive the following recursion 
\begin{equation}
\begin{array}{lll}
& ( \alpha\mathcal{L} \otimes \mathbf{I}_m  - 2\alpha^2\mathcal{L}^2 \otimes \mathbf{I}_m) (  X_{k+1}- X^{*} )\\& =  \mathcal{W}(X_k-X_{k+1})-\alpha  ( \mathcal{L}  \otimes \mathbf{I}_m) \big(  \Lambda_{k+2} - \Lambda^{*} \big)-Z_{k+1} \\&-\alpha \big(
\nabla  \widetilde{f} (X_k) -\nabla \widetilde{f} (X^{*})  \big)
- \alpha( \nabla \widetilde{f} (X^{*})+( \mathcal{L}  \otimes \mathbf{I}_m)\Lambda^{*} ) ,\nonumber
\end{array}
\end{equation}
or in the  alternative form
\begin{equation} 
\begin{array}{lll}
&   \mathcal{W}(X_{k+1}-X_k )+ \alpha  ( \mathcal{L}  \otimes I_m) \big(  \Lambda_{k+2} - \Lambda^{*} \big)\\&=
-( \alpha\mathcal{L} \otimes I_m  - 2\alpha^2\mathcal{L}^2 \otimes I_m) (  X_{k+1}- X^{*} )-Z_{k+1}\\&  \quad -\alpha \big(
\nabla  \widetilde{f} (X_k) -\nabla \widetilde{f} (X^{*})  \big)
- \alpha(\nabla  \widetilde{f} (X^{*})+( \mathcal{L}  \otimes I_m)\Lambda^{*} ) . \nonumber
\end{array}
\end{equation}
Then by \eqref{inequality01}   we derive
\begin{equation}\label{inequality1}
\begin{array}{lll}
&  \langle  \Lambda_{k+2}-\Lambda_{k+1},  \Lambda_{k+2}-\Lambda^{*}  \rangle
 +   \langle  X_{k+1}-X^*, \mathcal{W} (X_{k+1} -X_k)  \rangle
   \\
 &= - \langle  X_{k+1}-X^{*},   \big( ( \alpha\mathcal{L} - 2\alpha^2\mathcal{L}^2) \otimes \mathbf{I}_m\big) (  X_{k+1}- X^{*} ) \rangle
 \\&\quad-\alpha \langle  X_{k+1}-X^{*}, \nabla \widetilde{f} (X_k) - \nabla \widetilde{f} (X^{*})  \rangle  \\ &
\quad - \langle  X_{k+1}-X^{*}, Z_{k+1}+ \alpha( \nabla \widetilde{f} (X^{*}) +(\mathcal{L} \otimes  \mathbf{I}_m) \Lambda^{*} )  \rangle .
\end{array}
\end{equation}

By the definition of  the saddle point we have
 \begin{equation}\label{saddlepoint}
\phi(X^{*}, \Lambda) \leq \phi(X^{*}, \Lambda^{*})  \leq \phi(X , \Lambda^{*}) ~~ \forall  X \in \Omega, \Lambda \in \mathbb{R}^{mn}.
\end{equation}
Therefore,     $X^{*}$  minimizes $\widetilde{f} (X ) +  \langle X, ( \mathcal{L} \otimes \mathbf{I}_m) \Lambda^{*} 
\rangle$ over $\Omega$.
Since $\phi(X, \Lambda)$ is convex in $X \in \Omega$ for   each $\Lambda$,  by    noticing $X_k \in \Omega~~\forall k \geq 0$
 from  the optimal condition   \cite[ Proposition 1.1.8 ]{Bertsekas}  we derive
 \begin{equation}\label{inequality2}
\langle \nabla \widetilde{f} (X^{*}) +( \mathcal{L} \otimes \mathbf{I}_m) \Lambda^{*},  X_{k+1}-X^{*} \rangle  \geq 0 ~~ \forall k \geq 0.
\end{equation}

From    \eqref{primal}  \eqref{z} it follows that  $P_{\Omega} \{  X_{k+1}+ Z_{k+1}\}= X_{k+1}$,
and hence   $Z_{k+1} \in  N_{\Omega} \{ X_{k+1}\}$ by \eqref{normalcone}.
Then by  the definition of normal cone we obtain
 \begin{equation}\label{inequality3}
\langle  X_{k+1}-X^{*}, Z_{k+1}  \rangle  \geq 0.
\end{equation}

Then by  combining \eqref{inequality1} \eqref{inequality2} \eqref{inequality3} we derive
\begin{equation} 
\begin{array}{lll}
&  \langle  \Lambda_{k+2}-\Lambda_{k+1},  \Lambda_{k+2}-\Lambda^{*}  \rangle
 +   \langle  X_{k+1}-X_k, \mathcal{W} (X_{k+1} -X^{*})  \rangle    \\
 &\leq - \langle  X_{k+1}-X^{*},   \big( ( \alpha\mathcal{L} - 2\alpha^2\mathcal{L}^2) \otimes \mathbf{I}_m\big) (  X_{k+1}- X^{*} ) \rangle
\\& ~~~-\alpha \langle  X_{k+1}-X^{*}, \nabla \widetilde{f} (X_k) - \nabla \widetilde{f} (X^{*})  \rangle . \nonumber   \end{array}
\end{equation}
This incorporating with  \eqref{duct0}  yields
 \begin{equation}  \label{duct}
\begin{array}{lll}
&    V(X_{k+1}, \Lambda_{k+2})-V(X_k, \Lambda_{k+1})       \\
& \leq  -\| \Lambda_{k+2}-\Lambda_{k+1}\|^2 -   \langle X_{k+1}- X_k, \mathcal{W}  ( X_{k+1}- X_k)  \rangle \\&
~~ -2\alpha \langle  X_{k+1}-X^{*}, \nabla \widetilde{f} (X_k) - \nabla \widetilde{f} (X^{*})  \rangle \\&
~~ -2 \langle  X_{k+1}-X^{*},   \big( ( \alpha\mathcal{L} - 2\alpha^2\mathcal{L}^2) \otimes \mathbf{I}_m\big) (  X_{k+1}- X^{*} ) \rangle  . \end{array}
\end{equation}

Since $\mathcal{L}$ is  symmetric,
 there exists an  orthogonal matrix $\mathcal{U}$ such that $\mathcal{U}^T\mathcal{L} \mathcal{U}=diag\{0,\kappa_2,  \cdots, \kappa_n \}, $ and hence $\mathcal{U}^T\mathcal{L}^2 \mathcal{U}=diag\{0,\kappa_2^2, \cdots, \kappa_n^2 \}.  $
  Then by \eqref{wmatrix} we know that all  possible distinct  eigenvalues of    $\alpha\mathcal{L} - 2\alpha^2\mathcal{L}^2$ are $0, $ and $\alpha\kappa_i-2\alpha^2\kappa_i^2, i=2,\cdots, n$.  If  $ 0< \alpha \leq \frac{1}{2\kappa_n}  ,$
 then  $2\alpha\kappa_i \leq 1~~ \forall i=1,\cdots, n$,  and hence
  $ \alpha\kappa_i-2\alpha^2\kappa_i^2= \alpha\kappa_i(1- 2\alpha\kappa_i)\geq 0~ ~\forall i=1,\cdots, n.$
 Therefore,  for  any $\alpha$ with $0<\alpha \leq \frac{1}{2\kappa_n}$    the   matrix  $\alpha\mathcal{L} - 2\alpha^2\mathcal{L}^2$  is positive semi-definite, and   hence  by \eqref{duct} we derive
   \begin{equation}  \label{bound2}
\begin{array}{lll}
&    V(X_{k+1}, \Lambda_{k+2})-V(X_k, \Lambda_{k+1})      \\
& \leq  -\| \Lambda_{k+2}-\Lambda_{k+1}\|^2 -  \langle X_{k+1}- X_k, \mathcal{W}  ( X_{k+1}- X_k)  \rangle  \\&  ~~ ~-2\alpha \langle  X_{k+1}-X^{*}, \nabla \widetilde{f} (X_k) - \nabla \widetilde{f} (X^{*})  \rangle.  \end{array}
\end{equation}

 Let  the constant   $\alpha$ be such that  $ 0< \alpha \leq \frac{1}{2\kappa_n}  \textrm{~and ~}\alpha < \frac{3 }{2l_r}$.
In what follows,   we  show that  $   V(X_k, \Lambda_{k+1})    $  monotonously   decreases  and
$ d_k \leq r ~~\forall k \geq 0$  by induction.

We first  show  that $ d_1 \leq r $ and $  V(X_{1}, \Lambda_{2}) \leq V(X_0, \Lambda_{1})  $.
By the definition of the local  Lipschitz  constant  $ l_r$,  we know  that $\nabla \widetilde{f}(X)$ is  Lipschitz continuous on the compact  set $\{X \in \Omega:\| X-X^{*}\| \leq r \} $  with Lipschitz  constant $ l_r$.
Since    $\|X_0-X^{*}\| \leq r $  by  the definition of $ r $,    from Lemma \ref{lemma5} we see
 $$\langle X_{0}-X^{*}, \nabla  \widetilde{f} (X_0) - \nabla \widetilde{f} (X^{*}) \rangle  \geq  \frac{1}{ l_r} \| \nabla \widetilde{f} (X_0) - \nabla  \widetilde{f} (X^{*})  \|^2 .$$ 
 This incorporating with $xy \leq \frac{x^2}{4}+y^2$ leads to
\begin{equation}\label{inequality4}
\begin{array}{lll}
&  -\langle X_{1}-X^{*}, \nabla \widetilde{f} (X_0) - \nabla \widetilde{f} (X^{*})  \rangle  \\ &=
-\langle X_{0}-X^{*}, \nabla  \widetilde{f} (X_0) - \nabla \widetilde{f} (X^{*}) \rangle \\&\quad+  \langle -X_{1}+X_0, \nabla  \widetilde{f} (X_0) - \nabla  \widetilde{f} (X^{*})  \rangle \\
& \leq  -\frac{1}{ l_r} \| \nabla \widetilde{f} (X_0) - \nabla  \widetilde{f} (X^{*})  \|^2  + \frac{ l_r}{ 4}
 \| X_0-X_{ 1} \|^2\\&\quad+ \frac{1}{ l_r} \| \nabla \widetilde{f} (X_0) - \nabla   \widetilde{f} (X^{*})  \|^2
 \leq \frac{ l_r}{ 4} \| X_0-X_{1} \|^ 2. \end{array}
\end{equation}
Then from here by \eqref{bound2} we have
 \begin{equation}  \begin{array}{lll}
&V(X_{1}, \Lambda_{2}) - V(X_0, \Lambda_{1})   \\
& \leq  -\| \Lambda_{2}-\Lambda_{1}\|^2 -  \langle X_{1}- X_0, ( \mathcal{W} -\frac{\alpha  l_r}{2} I_{mn})( X_{1}- X_0) \rangle .
\end{array}
\nonumber \end{equation}

By \eqref{wmatrix} we know that
  all  possible  distinct  eigenvalues of  $\mathcal{W} -\frac{\alpha  l_r}{2 } I_{mn}$  are  $
1-\alpha\kappa_i+\alpha^2\kappa_i^2- \frac{\alpha  l_r}{2 },i=1,\cdots,n$.  Since   $  \alpha < \frac{3 }{2l_r},$
 we have $1-\alpha\kappa_i+\alpha^2\kappa_i^2- \frac{\alpha  l_r}{2 }
=(\frac{1}{2}-\alpha\kappa_i)^2+\frac{3}{4}- \frac{\alpha  l_r}{2 }>0 $. Thus,  $ \mathcal{W} -\frac{\alpha  l_r}{2} I_{mn}$ is positive definite. As a result,  $  V(X_{1}, \Lambda_{2}) \leq V(X_0, \Lambda_{1})  $. 
Then $  V(X_{1}, \Lambda_{2}) \leq  r^2 \lambda_{min} (\mathcal{M})$,  and hence $ d_1 \leq r $.

Assume    that $ d_p \leq r $ and   $   V(X_p, \Lambda_{p+1}) \leq    V(X_{p-1}, \Lambda_{p})    $     for $p =1,\cdots,k$.  Since \eqref{bound2} holds and $\| X_k-X^{*}\|\leq r$,  similar to the case  $k=0$,
we  can show  that     $   V(X_{k+1}, \Lambda_{k+2}) \leq    V(X_{k}, \Lambda_{k+1})    $ and $ d_{k+1} \leq r $.

In summary, by  the mathematical induction  we conclude that  $ d_k \leq r ~~\forall k \geq 0$, and  
 $   V(X_k, \Lambda_{k+1})    $    monotonously  decreases.
 
 Since  $\| X_k-X^{*}\| \leq r ~~\forall k \geq 0$,
  by the same procedure for deriving   \eqref{inequality4} we obtain     \begin{equation}
 -\langle X_{k+1}-X^{*}, \nabla \widetilde{f} (X_k) - \nabla \widetilde{f} (X^{*})  \rangle   \leq \frac{ l_r}{ 4} \| X_k-X_{k+1} \|^ 2~~\forall k \geq 0.  \nonumber
\end{equation}
Then  from here by  \eqref{bound2} we derive   
 \begin{equation} \label{deduce}
\begin{split}
& \quad  V(X_{k+1}, \Lambda_{k+2})-V(X_k, \Lambda_{k+1})  \leq  -\| \Lambda_{k+2}-\Lambda_{k+1}\|^2  \\&  - 
 \langle X_{k+1}- X_k,( \mathcal{W}  -\frac{\alpha  l_r}{2 } I_{mn} ) ( X_{k+1}- X_k)  \rangle \leq 0   .
\end{split} \end{equation}
 Thus, we conclude that    $   V(X_k, \Lambda_{k+1})    $   converges since
it is  nonnegative.    Summing  up   both sides of   \eqref{deduce}  from $0$ to $p$ we derive
 \begin{equation}\label{deduce01}
 \begin{array}{lll}
&    V(X_{p+1}, \Lambda_{p+2}) -   V(X_0, \Lambda_{1})         \leq - \sum_{k=0}^p  \| \Lambda_{k+2}-\Lambda_{k+1}\|^2 \\& ~~~~- \sum_{k=0}^p   \langle X_{k+1}- X_k,( \mathcal{W}  -\frac{\alpha  l_r}{2 } I_{mn} ) ( X_{k+1}- X_k)  \rangle  .
\end{array}
\end{equation}
Then by letting $p \rightarrow \infty$ we have
\begin{align}
& \sum_{k=0}^{\infty} \langle X_{k+1}- X_k,( \mathcal{W}  -\frac{\alpha  l_r}{2 } I_{mn} ) ( X_{k+1}- X_k)  \rangle     < \infty, 
\label{sum01}    \\
&\textrm{and ~ } \sum_{k=1}^{\infty}   \| \Lambda_{k+1}- \Lambda_k \| ^2  < \infty.\label{sum02}
\end{align}

Consequently, we   derive
$\lim\limits_{k \rightarrow  \infty} (X_{k+1}-X_k)=0$ by \eqref{sum01} since   $\mathcal{W} -\frac{\alpha  l_r}{2} I_{mn}$   is positive definite, and
 $\lim\limits_{k \rightarrow  \infty} (\Lambda_{k+1} -\Lambda_k)=0$ by  \eqref{sum02}.  By  convergence of    $   V(X_k, \Lambda_{k+1})    $    we conclude that $X_k$ and $  \Lambda_k$ contain    convergent subsequences
  $\{X_{ n_k}\}$  and $\{ \Lambda_{n_k}\}$ to some   limits   $X^0$  and $\Lambda^0$, respectively.
   Since $\lim\limits_{k \rightarrow  \infty} (X_{n_k+1}-X_{n_k})=0$ and
  $\lim\limits_{k \rightarrow  \infty} (\Lambda_{n_k+1}-\Lambda_{n_k})=0$,  by noticing that  $P_{\Omega}\{x\} $ is continuous in $x$ and
 $\lim\limits_{k \rightarrow  \infty}  X_{n_k} =X^0,~\lim\limits_{k \rightarrow  \infty}  \Lambda_{n_k} =\Lambda^0$,
   from \eqref{primal} \eqref{dual} we derive
 \begin{align}
& X^0= P_{\Omega} \{  X^0 -  \alpha \nabla  \widetilde{f} (X^0)  -\alpha ( \mathcal{L}   \otimes \mathbf{I}_m)  \big( \Lambda^0 +  X^0 \big)  \}, \label{primal1} \\
&  ( \mathcal{L}  \otimes \mathbf{I}_m)  X^0  =0. \label{dual1}
\end{align}

Then  from \eqref{primal1} \eqref{dual1} by \eqref{normalcone} we see $  \alpha  \big ( \nabla \widetilde{f} (X^0) +(\mathcal{L}   \otimes \mathbf{I}_m) \Lambda^0  \big) = \alpha  \big ( \nabla \widetilde{f} (X^0) +(\mathcal{L}   \otimes \mathbf{I}_m)( \Lambda^0 +  X^0 ) \big)  \in  - N_{\Omega}(X^0),$
and hence  by the definition of normal cone we conclude
$$ \langle \nabla \widetilde{f} (X^0) +(\mathcal{L}   \otimes \mathbf{I}_m)  \Lambda^0 , X-X^0 \rangle \geq 0~~ \forall X \in \Omega .$$
Since    $\phi(X,\Lambda)= \widetilde{f}(X)+\langle\Lambda, ( \mathcal{\bar{L}} \otimes \mathbf{I}_m ) X\rangle$ is convex
in $X \in \Omega$ for  each $\Lambda \in \mathbb{R}^{mn} $,  by \eqref{e1} we have
\begin{equation}
\begin{array}{lll}
\phi(X,\Lambda^0)& \geq \phi(X^0,\Lambda^0)+   \langle \nabla \widetilde{f} (X^0) +(\mathcal{L}   \otimes \mathbf{I}_m ) \Lambda^0   , X-X^0 \rangle  \\&\geq \phi(X^0,\Lambda^0)  ~~ \forall X \in \Omega. \nonumber
\end{array}
\end{equation}
From \eqref{dual1} we see   $\phi(X^0,\Lambda^0) =\phi(X^0,\Lambda)=\widetilde{f}(X^0)  ~~\forall \Lambda \in \mathbb{R}^{mn}$, and hence   by definition we know  $(X^0, \Lambda^0)$ is a  saddle point of   $\phi(X,\Lambda)$ in $ \Omega \times \mathbb{R}^{mn}$.
  Thus,  by Lemma \ref{lem3} we see that $X^0$ is an optimal solution  to the problem \eqref{problem2}.

  Since $ \| \Lambda_{k+1}-\Lambda^{*} \|^2 +(X_{k} -X^{*})^T  \mathcal{W} (X_{k} -X^{*}) $  converges,
  by setting $(X^{*}, \Lambda^{*})=(X^0, \Lambda^0)$ from $
\lim\limits_{k \rightarrow  \infty}  X_{n_k}=X^0 $ and $ \lim\limits_{k \rightarrow  \infty} \Lambda_{n_k}=\Lambda^0$
    we conclude that  $ \| \Lambda_{k+1}-\Lambda^0 \|^2 +(X_{k} -X^0)^T  \mathcal{W} (X_{k} -X^0) $ converges to zero.
Therefore,
$$ \lim_{k \rightarrow \infty}  X_k= X^0, ~~ \lim_{k \rightarrow \infty}  \Lambda_k= \Lambda^0.$$
Thus,   by Remark \ref{r1}   we conclude that  $ \lim\limits_{k \rightarrow \infty} x_{i,k}=\lim\limits_{k \rightarrow \infty} x_{j,k}=x^* ~~\forall i,j \in \mathcal{V}$ for some $x^* \in \Omega_o^*.$  
\hfill $\blacksquare$

 \subsection{Proof of Theorem \ref{thm2}}
\textbf{Proof:} (i) Summing up both sides of \eqref{dual}  from $0$ to $p$   leads to
  \begin{equation} \Lambda_{p+1}-\Lambda_0= \sum_{k=0}^p \alpha   ( \mathcal{L}  \otimes \mathbf{I}_m)  X_k =(p+1)
 \alpha   ( \mathcal{L}  \otimes \mathbf{I}_m)  \bar{X}_p , \nonumber
\end{equation}
and hence  \eqref{con} holds.

(ii) When $\Omega_i =\mathbb{R}^m~\forall i\in \mathcal{V},$  the equation   \eqref{primal} turns   to
 \begin{equation}\label{primal2}
 X_{k+1}=  X_k -  \alpha     \nabla \widetilde{f} (X_k)   -\alpha ( \mathcal{L}   \otimes \mathbf{I}_m)  \big( \Lambda_k +  X_k \big) .
\end{equation}
By  \eqref{deduce01} we see
  \begin{equation}\label{deduce02}
  \begin{array}{lll}
& \sum_{k=0}^p \big ( \| \Lambda_{k+2}-\Lambda_{k+1}\|^2 +   \langle X_{k+1}- X_k,( \mathcal{W}  -\frac{\alpha  l_r}{2 } I_{mn} ) ( X_{k+1}- X_k)  \rangle   \\&  \leq - V(X_{p+1}, \Lambda_{p+2}) +   V(X_0, \Lambda_{1})     . \nonumber
\end{array}
\end{equation}
Then by noticing that $\mathcal{W} -\frac{\alpha  l_r}{2} I_{mn}$   is positive definite we derive 
\begin{align} 
& \sum_{k=0}^p   \| \Lambda_{k+2}-\Lambda_{k+1}\|^2    \leq - V(X_{p+1}, \Lambda_{p+2}) +   V(X_0, \Lambda_{1})    \label{bound0} , \\&
\textrm{and  }\sum_{k=0}^p   \| X_{k+1}- X_k\|^2  \leq  \frac{ - V(X_{p+1}, \Lambda_{p+2}) +   V(X_0, \Lambda_{1})  }{\lambda_{min}( \mathcal{W} -\frac{\alpha  l_r}{2 } I_{mn } )} . \label{bound1}
\end{align}  

 By \eqref{primal2}   we have    \begin{equation}\label{eprimal2}
\begin{array}{lll}
&\|  X_{k+1} - X^{*}\|^2    = \|X_k-X^{*}\|^2+ \|X_{k+1}-X_k\|^2 \\&+2(   X_k-X^{*})^T (X_{k+1}-X_k)
= \|  X_k -X^{*}\|^2  \\&+2 \alpha  
\langle  X^{*} -X_k,     \nabla \widetilde{f} (X_k)   + ( \mathcal{L}   \otimes \mathbf{I}_m)  \big( \Lambda_k +  X_k \big)\rangle 
\\& + \|X_{k+1}-X_k\|^2.
\end{array}
\end{equation}
Noticing   $\phi(X, \Lambda)$ is convex in $X \in \Omega$ for   any $\Lambda$,  by \eqref{e1} we   have
 \begin{equation}\label{estimate}
\langle  X^{*} -X_k,   \nabla \widetilde{f} (X_k)   + ( \mathcal{L}   \otimes \mathbf{I}_m)   \Lambda_k  \rangle \leq
\phi(X^{*}, \Lambda_k)-\phi(X_k, \Lambda_k) .
\end{equation}
Since  $\Lambda_k$ is bounded and $X^{*}$ is  an optimal solution to the problem \eqref{problem2},
   by  \eqref{consensus} we see $\phi(X^{*}, \Lambda_k) =\widetilde{f}(X^{*})=f^{*}. $
Since  $\mathcal{L}$  is positive semi-definite by Lemma \ref{lem1},  from   \eqref{consensus}  it follows that
 $$ \langle X^{*} -X_k,  ( \mathcal{L}   \otimes \mathbf{I}_m)   X_k \rangle =
  - \langle X_k ,   ( \mathcal{L}   \otimes \mathbf{I}_m)   X_k \rangle  \leq 0.$$
  Thus,  from here by   \eqref{eprimal2} \eqref{estimate} we  derive
  \begin{equation}
 \|  X_{k+1} - X^{*}\|^2   \leq  \|  X_k -X^{*}\|^2 + 2\alpha(  f^{*}-\phi(X_k, \Lambda_k) ) + \|X_{k+1}-X_k\|^2, \nonumber
\end{equation}
and hence
  \begin{equation} \phi(X_k, \Lambda_k) - f^{*} \leq  \frac{1}{ 2\alpha}(\|  X_k -X^{*}\|^2  -\|  X_{k+1} - X^{*}\|^2   + \|X_{k+1}-X_k\|^2).\nonumber
\end{equation}
Summing up both sides of this  inequality   from $0$ to $p$ for $p  \geq 1$ we obtain
 \begin{equation} 
\begin{array}{lll}
&\sum_{k=0}^p ( \phi(X_k, \Lambda_k) - f^{*} ) \\& \leq  \frac{1}{ 2 \alpha} \big (\|  X_0 -X^{*}\|^2  -\|  X_{p+1} - X^{*}\|^2 + \sum_{k=0}^p \|X_{k+1}-X_k\|^2\big). \nonumber
\end{array}
\end{equation}
From here by the  convexity of $\widetilde{f}(X)$ we derive   \begin{equation} \label{1sum0}
\begin{array}{lll}
 \widetilde{f}(\bar{X}_p) &\leq\frac{1}{p+1} \sum_{k=0}^p   \widetilde{f}(X_k)
\leq \frac{1}{p+1}   \sum_{k=0}^p \phi(X_k, \Lambda_k) \\& - \frac{1}{p+1}   \sum_{k=0}^p X_k^T (  \mathcal{L} \otimes \mathbf{I}_m)  \Lambda_k \\&   \leq f^{* } + \frac{1}{2\alpha (p+1)}  \big(  \|  X_0 -X^{*}\|^2  -\|  X_{p+1} - X^{*}\|^2 \\&+ \sum_{k=0}^p \|X_{k+1}-X_k\|^2  \big) \\&  - \frac{1}{p+1}   \sum_{k=0}^p \langle X_k,(  \mathcal{L} \otimes \mathbf{I}_m)  \Lambda_k  \rangle .\end{array}
\end{equation}

 We  now give an upper bound for  $- \frac{1}{p+1}   \sum_{k=0}^p \langle X_k,(  \mathcal{L} \otimes \mathbf{I}_m)  \Lambda_k  \rangle$.
 By \eqref{dual} we have
  \begin{equation}
\label{ }
\| \Lambda_{k+1}\|^2 = \| \Lambda_k\|^2 + 2 \alpha  \langle  \Lambda_k,  ( \mathcal{L}  \otimes \mathbf{I}_m)  X_k \rangle
+\| \alpha   ( \mathcal{L}  \otimes \mathbf{I}_m)  X_k  \|^2. \nonumber
\end{equation}
Thus,     $-\langle X_k,(  \mathcal{L} \otimes \mathbf{I}_m)  \Lambda_k  \rangle
=\frac{1}{ 2\alpha} ( \| \Lambda_k\|^2 -\| \Lambda_{k+1}\|^2
+\| \Lambda_{k+1}-\Lambda_k \|^2)$,  and hence
\begin{equation}\label{1sum1}
\begin{split}
& -\frac{1}{p+1}   \sum_{k=0}^p  \langle X_k,(  \mathcal{L} \otimes \mathbf{I}_m)  \Lambda_k  \rangle \\&=\frac{1}{2\alpha(p+1)} ( \| \Lambda_0\|^2 -\| \Lambda_{p+1}\|^2
+ \sum_{k=0}^p \| \Lambda_{k+1}-\Lambda_k \|^2).
\end{split}
\end{equation}
By substituting \eqref{1sum1}  into  \eqref{1sum0}    we derive
\begin{equation} \begin{array}{lll}
 & \widetilde{f}(\bar{X}_p) \leq f^{* } \\&+ \frac{1}{2\alpha (p+1)}  \Big(  \|  X_0 -X^{*}\|^2  -\|  X_{p+1} - X^{*}\|^2  +\| \Lambda_0\|^2 -\| \Lambda_{p+1}\|^2   \Big ) \\&
 +  \frac{1}{2\alpha (p+1)}  \Big(   \sum_{k=0}^p  \|X_{k+1}-X_k\|^2   +  \sum_{k=0}^p  \| \Lambda_{k+1}-\Lambda_k\|^2  \Big) . \nonumber\end{array}
\end{equation}
Then from here by  \eqref{bound0}   \eqref{bound1}  we obtain   \eqref{sum0}.

(iii) By  \eqref{consensus} \eqref{saddlepoint} we derive $\phi( \bar{X}_p ,\Lambda^{*}) \geq\phi( X^{*} ,\Lambda^{*})=\widetilde{f}(X^{*}) =f^{*}$,
and hence  for  any dual solution $\Lambda^{*}$
\begin{equation} 
\begin{split}
& \widetilde{f}(\bar{X}_p) =\widetilde{f}(\bar{X}_p)  + \langle \bar{X}_p , (\mathcal{L} \otimes \mathbf{I}_m) \Lambda^{*} \rangle
- \langle \bar{X}_p,  (\mathcal{L} \otimes \mathbf{I}_m ) \Lambda^{*} \rangle \\& =\phi( \bar{X}_p ,\Lambda^{*})-\langle \bar{X}_p,  (\mathcal{L} \otimes \mathbf{I}_m ) \Lambda^{*} \rangle  \geq   f^{*}-\langle \bar{X}_p,  (\mathcal{L} \otimes \mathbf{I}_m ) \Lambda^{*} \rangle . \nonumber \end{split}
\end{equation} Then by \eqref{con} we derive  \eqref{sum00}.
\hfill $\blacksquare$

 \section{Numerical Simulations}
 In this section, we give two  numerical examples to demonstrate the   obtained  theoretic results.

\begin{figure}     \centering
  \includegraphics[width=3in,height=3.3in]{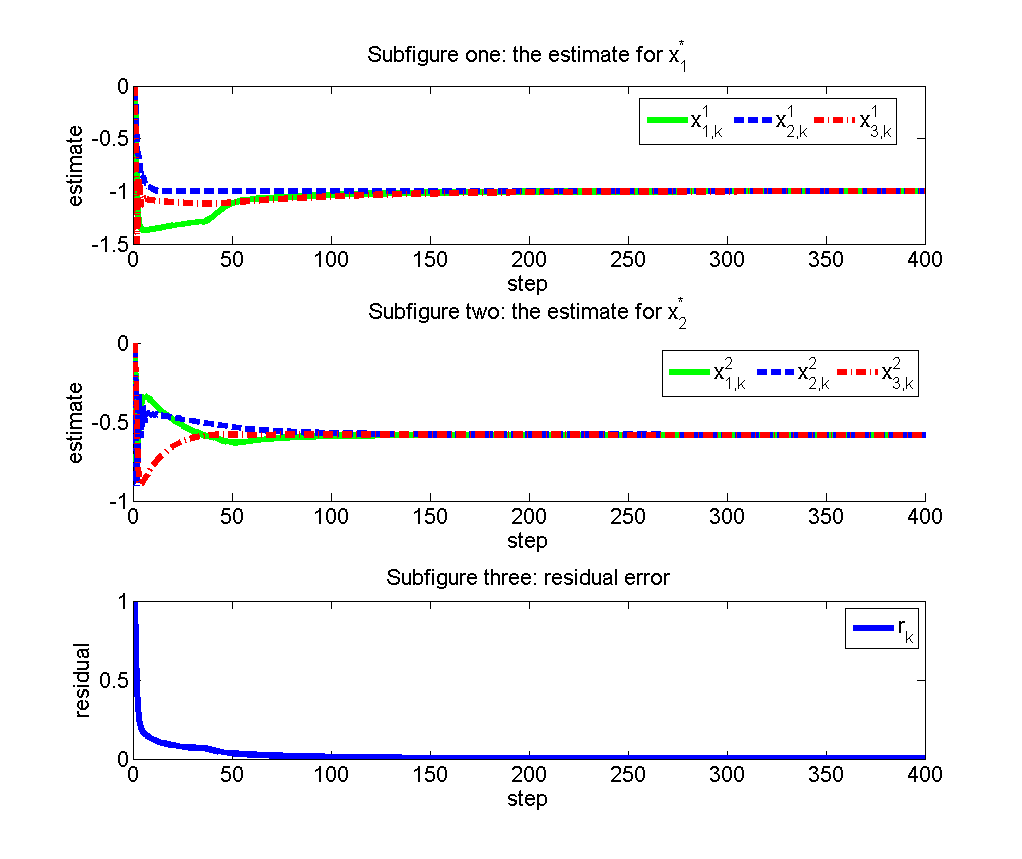}
   \caption{ The estimates and the residual } \label{A3}
\end{figure}

\textbf{ Example 5.1.} This example shows  that    the  primal-dual algorithm with constant step size can produce the accurate  estimates  for     the constrained  optimization problem,
  where   gradients of the cost functions are only locally Lipschitz continuous,   and  the  agents are equipped with
  different   constraint sets.  Besides, some of the constraint sets are not compact.

 Consider  an  undirected  network of three agents with   edge set
 $\mathcal{E}_{\mathcal{G}}=\{ (1,3),(2,3),(1,1),(2,2),(3,3)\}$.
 Objective functions for the  agents  are as follows:
 \begin{equation}
\begin{split}
&f_1(x_1,x_2)=\frac{x_1^2}{2}+3x_1+x_2^2+2x_2+x_1x_2+0.5e^{x_1+x_2}, \\
&f_2(x_1,x_2)=x_1^2+2x_1+2x_2^2+2x_2+x_1x_2+e^{x_2}, \\
&f_3(x_1,x_2)=2x_1^2+4x_1+x_2^2+2x_2+e^{x_1},
\end{split}
\end{equation}
while the constraint sets for agents are $\Omega_1=\{(x_1,x_2) \in \mathbb{R}^2: x_1^2+x_2^2\leq 2\},
~\Omega_2=\{(x_1,x_2)\in \mathbb{R}^2: x_1 \geq -1\},$ and $\Omega_3=\{(x_1,x_2)\in \mathbb{R}^2:   x_2 \leq -0.5\}.$
Denote  the optimal solution by   $(x_1^{*},x_2^{*}) $, which is at the boundary of the global constraint set.

Let $\{ x_{i,k}\}$ and  $\{ \lambda_{i,k}\}$   be produced by the algorithm \eqref{fixed}
 with initial values $x_{i,0}=\mathbf{0},~ \lambda_{i,0}=\mathbf{0},~i=1,2,3$, and  $\alpha=0.4$.
 Denote by $x_{i,k}^1$ and $x_{i,k}^2$ the estimates for $ x_1^{*}$ and $ x_2^{*} $ by agent $i$ at time $k,$ respectively.  Note that the primal-dual solution pair $(X^{*}, \Lambda^{*})$  satisfies \eqref{primal1} and \eqref{dual1}.
Define the residual of the optimal condition as
$r_k=col\{ X_{k+1}-X_k,  ( \mathcal{L}  \otimes \mathbf{I}_m)  X_k\}. $ The local  estimates of
all agents  and 2-norm of  the residual $r_k$
are shown in  Figure  \ref{A3}. From the figure    it is seen  that the estimates for all agents   converge to the same  optimal solution.

\begin{figure}
    \centering
  \includegraphics[width=3in]{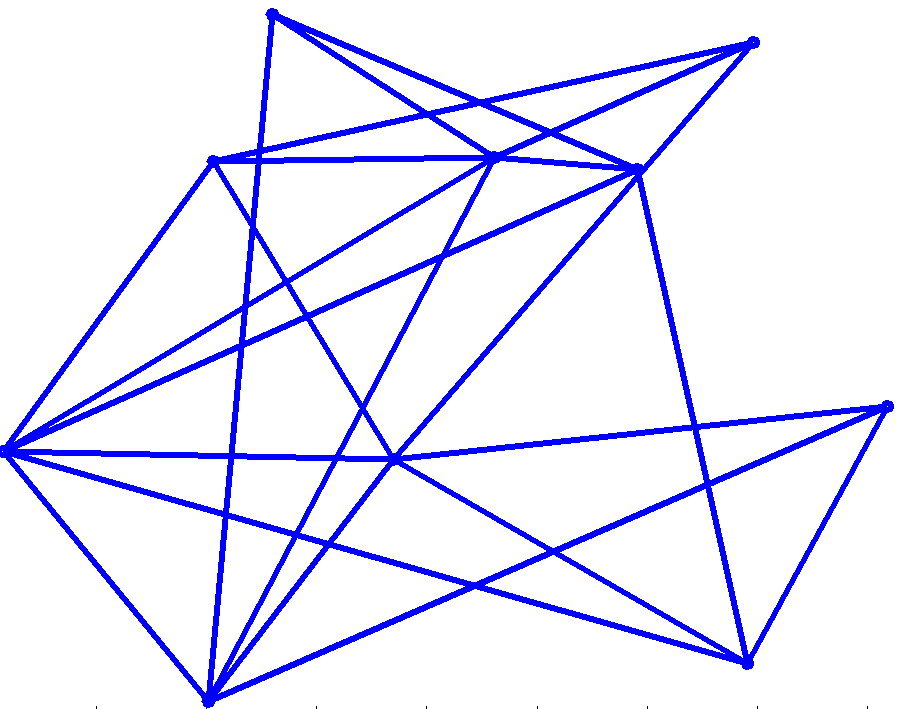}
   \caption{The communication topology} \label{Figure A1}
\end{figure}

\textbf{Example 5.2. }Consider a randomly generated undirected network of  $n=10$ agents,  where  each agent has  an   average degree 4.  Each agent $i \in \mathcal{V}$ is assigned with a huber loss function $f_i:\mathbb{R}\rightarrow \mathbb{R}$ with
  $$f_i(x)= \begin{cases} & \frac{1}{2}(x-a_i)^2 ,~\textrm{if}~ | x-a_i| \leq 1,\\ &   |x-a_i|-\frac{1}{2} , \textrm{  otherwise.}   
  \end{cases} $$
   For any $i \in \mathcal{V},$ $a_i$  is generated  according to  the uniform distribution over the interval  $[1.5,2.5]$. The optimal solution  of $f(\cdot)=\sum_{i=1}^n f_i(\cdot)$   is denoted by $x^{*}$. We compare  the  primal-dual algorithm  \eqref{fixed} with the existing ones  by    this  example.

  Set $\alpha=0.8.$ Denote by $x_{i,k}$   the estimate  for $x^{*}$  given by agent $i$ at time $k$  with  the initial values
 $x_{i,0}=\mathbf{0}~~\forall i \in \mathcal{V}$.
The simulation is for the case where the  communication topology   is  shown in Figure \ref{Figure A1},
  and   the entries of the  adjacency  matrix   $\mathcal{A_{G}}$  are Metropolis wights \cite{Metropolis}.  With this  $\mathcal{A_G}$ we carry out the   simulations  for the  primal-dual  algorithm \eqref{fixed} with
  $\lambda_{i,0}=\mathbf{0}~\forall i \in \mathcal{V}$,
    for the   DGD algorithm  \cite{optimization0},  for  EXTRA    \cite{Ling} with   constant  step size $\alpha$,
  and for  the distributed Nesterov gradient (D-NG)  algorithm    in \cite{Moura}.
 The DGD algorithm runs separately for   three cases:   constant step size $\alpha$, diminishing step sizes $\alpha_k=\alpha/k^{0.75}$,  and $\alpha_k=\alpha/k^{0.4}$. The  D-NG algorithm, i.e.,  equations (2)-(4) in \cite{Moura},  is run with $c=\alpha$ and $y_{i,0}=\mathbf{0}~\forall i  \in \mathcal{V}.$

\begin{figure}
      \centering
       \includegraphics[width=3in]{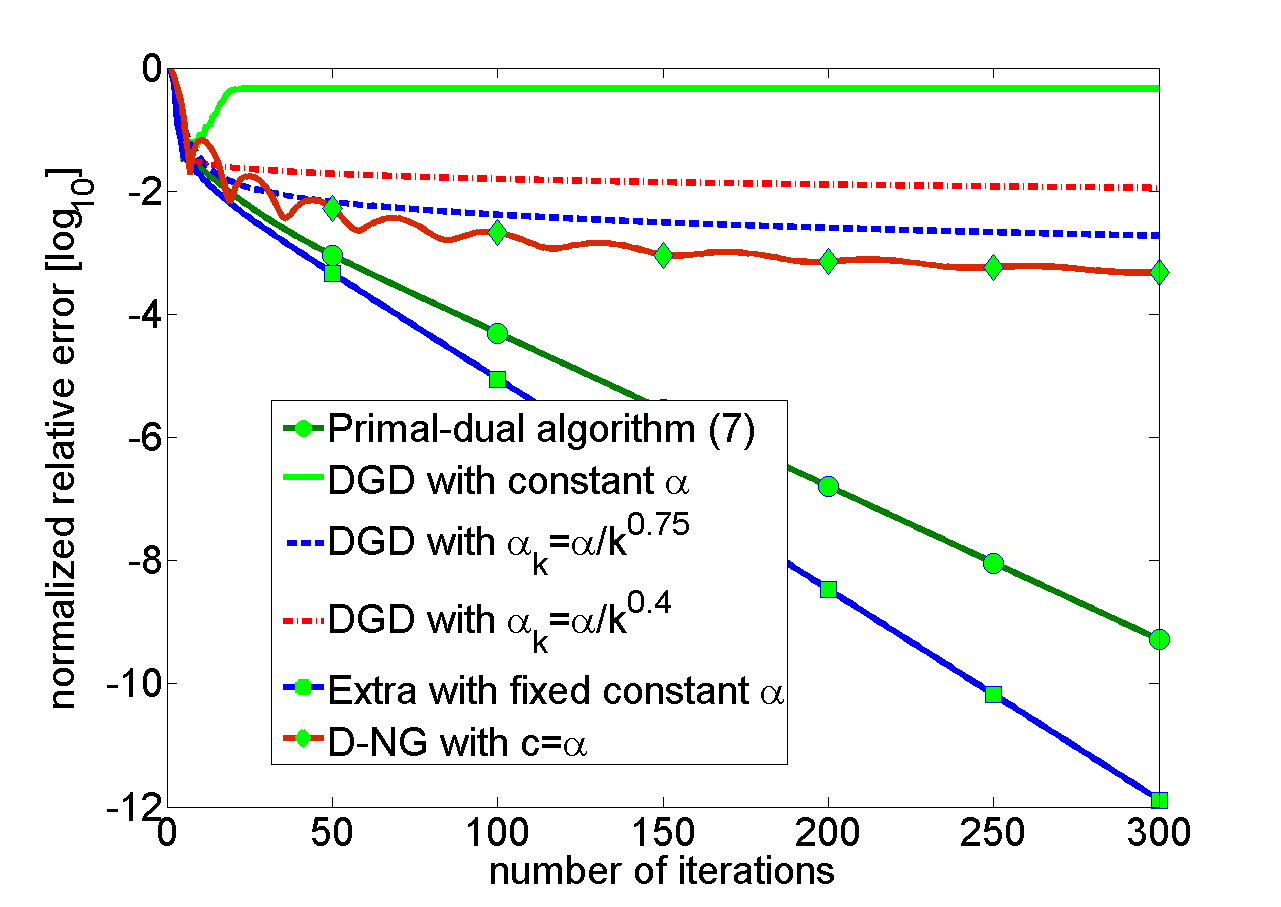}
 \caption{The normalized relative error} \label{Figure A2}
\end{figure}

Denote  by
 $e_k=\frac{\|X_k-\mathbf{1} \otimes x^{*}\|}{\|X_0-\mathbf{1} \otimes x^{*}\|} $    the normalized relative error, where  $X_k=col\{x_{1,k},\cdots, x_{n,k}\}$.
  The numerical results  are shown in Figure \ref{Figure A2},
where  the horizontal axis  denotes the number of iterations $k$ and  the vertical axis denotes $\textrm{log}_{10} (e_k)$.
From the figure  it is seen  that  the DGD  algorithms with decreasing  step sizes converge to the optimal solution but the
rate of convergence are the   slowest  in comparisons with other methods. 
It is also seen that DGD   with constant step size quickly   approaches to the neighborhood of the optimal    solution.
The estimates generated by  D-NG \cite{Moura},  by   the algorithm \eqref{fixed}, and by EXTRA \cite{Ling}  all converge  to the optimal solution. Besides,  the  algorithm \eqref{fixed}     brings  a  satisfactory convergence rate for the unconstrained problem as well.

 \section{Conclusion}
In the paper,   a  distributed primal-dual  algorithm is proposed  for  multiple agents in a
 network to    minimize the sum
of individual cost functions subject to a global constraint, which is the intersection of the  local constraints.
   The proposed algorithm with constant step size makes  the  estimates   of all agents converge to  the same  optimal solution  and  achieve  the  convergence  rate $O(\frac{1}{k})$ when there is no constraint.
 The effectiveness  and the priority of the proposed algorithm have  been demonstrated by two numerical examples.

For further research, it is of interest to consider   the primal-dual algorithm for  stochastic optimization, and  to see
if some desired properties taking place for   the deterministic   still remain true.

   \end{document}